\documentclass[11pt]{article}
\usepackage{amsmath,amssymb,theorem}
\usepackage{graphicx}
\usepackage{caption}

\usepackage{subfigure}
\usepackage{enumerate}

\usepackage{psfrag}
\usepackage{mathrsfs}
\usepackage{tikz-cd}
\usepackage{bbm}

\usepackage{epstopdf}
\usepackage{float}

\theorembodyfont{\rmfamily}
\newtheorem{thm}{Theorem}[section]
\newtheorem{conj}[thm]{Conjecture}
\newtheorem{lem}[thm]{Lemma}%[section]
\newtheorem{prop}[thm]{Proposition}%[section]
\theorembodyfont{\rmfamily}
\def\pf{\bigskip\noindent {\bf Proof.}~~}

\def\mytextindent#1{\indent\llap{#1\enspace}\ignorespaces}
\def\fattextindent#1{\indent\indent\llap{#1\enspace}\ignorespaces}
\def\myitem{\par\hangindent\parindent\mytextindent}
\def\myitemitem{\par\hangindent\parindent\fattextindent}
\def\proofsquare{  \bs\hfill $\blacksquare$}

\def\bs{\bigskip}
\def\less{\backslash}
\def\dfn#1{{\sl #1}}

\begin{document}
\title{Clique minors in double-critical  graphs }
\author{Martin Rolek\thanks {Current address: Department of Mathematics, College of William and Mary.  E-mail address: msrolek@wm.edu.}
~and~ Zi-Xia Song\thanks{Corresponding author. E-mail address: Zixia.Song@ucf.edu.}\\
Department of Mathematics\\
University of Central Florida\\
 Orlando, FL 32816
}
%\date{March 19, 2016}
\maketitle
\begin{abstract}
A connected $t$-chromatic graph $G$ is  \dfn{double-critical} if   $G -\{u , v\}$ is $(t-2)$-colorable  for each edge $uv\in E(G)$.   A long-standing  conjecture of Erd\H{o}s and Lov\'asz that the complete graphs are the only double-critical $t$-chromatic graphs  remains open for all $t\ge6$. Given the difficulty in settling Erd\H{o}s and Lov\'asz's conjecture and motivated by the well-known Hadwiger's conjecture, Kawarabayashi, Pedersen and Toft  proposed a weaker conjecture that every double-critical $t$-chromatic graph contains a $K_t$  minor and  verified their conjecture  for $t\le7$.   Albar and Gon\c calves recently proved that  every double-critical $8$-chromatic graph contains a $K_8$  minor, and their proof is computer-assisted. In this paper   we prove that every  double-critical $t$-chromatic graph contains a $K_t$  minor for all $t\le9$. Our proof for $t\le8$ is shorter and  computer-free.  
\end{abstract}

\baselineskip=18pt
\section{Introduction}

All graphs in this paper are finite and simple.  
For a graph $G$ we use $|G|$, $e(G)$, $\delta (G)$ to denote the number
of vertices, number of edges and minimum degree of $G$, respectively.
The degree of a vertex $v$ in a graph is denoted by $d_G(v)$ or
simply $d(v)$.   For a subset $S$ of $V(G)$,
the subgraph induced by $S$ is denoted by
$G[S]$ and $G - S =G[V(G)\setminus S]$.
%For a subgraph $H$ of $G$, $G - V(H) =G[V(G)\setminus V(H)]$.  
  If $G$ is a graph and $K$ is a subgraph of $G$, then by $N(K)$ we denote
the set of vertices of $V(G)\setminus V(K)$ that are adjacent to a vertex of $K$.
If $V(K)=\{x\}$, then we use $N(x)$ to denote $N(K)$.  By abusing
notation we will also denote by $N(x)$ the graph induced by the set
$N(x)$.  We define $N[x]=N(x)\cup \{x\}$, and similarly will use the
same symbol for the graph induced by that set.  %If $x,y$ are adjacent vertices of a graph $G$, then we denote by $G/xy$ the graph obtained from $G$ by contracting the edge $xy$ and deleting all resulting parallel edges.  
If $u,v$ are distinct nonadjacent vertices of a graph $G$, then by
$G+uv$ we denote the graph obtained from $G$ by adding an edge
with ends $u$ and $v$.  If $u,v$ are adjacent or equal, then we define
$G+uv$ to be $G$.

A graph  $H$ is a \dfn{minor} of a graph $G$ if  $H$ can be
 obtained from a subgraph of $G$ by contracting edges.  We write $G\ge H$ if 
$H$ is a minor of $G$.
In those circumstances we also say that  $G$ has an $H$ \dfn{minor}.  
A connected graph $G$ is called \dfn{double-critical} if  for any edge $uv \in E(G)$, we have $\chi(G -\{u , v\}) = \chi(G) - 2$.   The following long-standing  \dfn{Double-Critical Graph Conjecture} is due to  Erd\H{o}s and Lov\'asz~\cite{Erdos1968}.
\begin{conj}\label{dc}{\bf Double-Critical Graph Conjecture} (Erd\H os and Lov\'asz~\cite{Erdos1968}) \, 
For every integer $t\ge1$, the only double-critical $t$-chromatic graph is $K_t$.
\end{conj}  

\noindent Conjecture~\ref{dc}   is a special case of the so-called Erd\H{o}s-Lov\'asz Tihany Conjecture~\cite{Erdos1968}. It is 
 trivially true for $t\le3$ and reasonably easy for  $t = 4$.   Mozhan \cite{Mozhan1987} and Stiebitz \cite{Stiebitz1987} independently proved Conjecture~\ref{dc}  for $t=5$.

\begin{thm}\label{dc5c} (Mozhan~\cite{Mozhan1987}; Stiebitz~\cite{Stiebitz1987})
The only double-critical $5$-chromatic graph is $K_5$.
\end{thm}

 Conjecture~\ref{dc} remains open for all $t \ge 6$.  Given the difficulty in settling Conjecture~\ref{dc}  and motivated by the well-known Hadwiger's conjecture~\cite{hc},  Kawarabayashi, Pedersen and Toft  proposed  a weaker conjecture.

\begin{conj}\label{dcminor} (Kawarabayashi, Pedersen and Toft~\cite{KPT2010})
For every integer $t\ge1$, every double-critical $t$-chromatic graph contains a $K_t$ minor.
\end{conj}

\noindent Conjecture~\ref{dcminor}  is a weaker version of Hadwiger's conjecture~\cite{hc}, which states that for every integer $t \ge 1$, every $t$-chromatic graph contains a $K_t$ minor.
Conjecture~\ref{dcminor}  is true  for $t \le 5$ by Theorem~\ref{dc5c}.   In the same paper~\cite{KPT2010}, Kawarabayashi, Pedersen and Toft verified  their conjecture for $t \in\{ 6, 7\}$. 

\begin{thm}\label{dc7c} (Kawarabayashi, Pedersen and Toft~\cite{KPT2010})
For every integer $t\le 7$, every double-critical $t$-chromatic graph contains a $K_t$ minor.
\end{thm}

 Recently,  Albar and Gon\c calves~\cite{AG2015} announced a proof for the case $t=8$.  

\begin{thm}\label{dc8c} (Albar and Gon\c calves~\cite{AG2015})
Every double-critical $8$-chromatic graph has a  $K_8$ minor.
\end{thm}

  Our main result is the following next step. 

\begin{thm}\label{dc9c}
For integers $k, t$ with $1\le k\le 9$ and $t\ge k$,  every   double-critical $t$-chromatic  graph  contains a $K_k$  minor. 
\end{thm}

We actually prove a much stronger result, the following.
\begin{thm}\label{strongk89}
For $k \in \{6, 7, 8, 9\}$, let $G$ be a  $(k-3)$-connected graph  with  $k+1\le\delta(G)\le 2k-5$.   If  every edge of $G$ is contained in at least $k-2$ triangles and for any    
  minimal separating set $S$ of $G$ and any $x\in S$, $G[S\less \{x\}]$   is not  a clique, 
then $G\ge K_k$. 
\end{thm}

  Theorem~\ref{dc9c} follows directly from Proposition~\ref{prop} (see below) and  Theorem~\ref{strongk89}.
Our proof of Theorem~\ref{strongk89}  closely follows the proof of the extremal function for $K_9$ minors by Song  and Thomas~\cite{SongThomas2006} (see Theorem~\ref{k9} below).
Note that the proof of Theorem~\ref{dc7c} for $k=7$ is about ten pages long and the proof of Theorem~\ref{dc8c} is  computer-assisted.
Our proof of Theorem~\ref{dc9c} is much shorter and computer-free for $k\le8$.
For $k=9$, our proof is computer-assisted as  it applies a computer-assisted lemma from~\cite{SongThomas2006} (see Lemma~\ref{k7k1} below).
Note that a computer-assisted proof of Theorem~\ref{strongk89}  for all  $k\le8$ (and hence computer-assisted proofs of Theorem~\ref{dc7c} and Theorem~\ref{dc8c}) follows directly from  Theorem~\ref{strongk89} for $k=9$.
(To see that, let $G$ and $k\le 8$ be as in  Theorem~\ref{strongk89}, and let $H$ be obtained from $G$ by adding $9-k$ vertices, each adjacent to every other vertex of the graph.
Then $H$ is $6$-connected and satisfies all the  other conditions as stated in Theorem~\ref{strongk89}.
Thus $H\ge K_9$ and so $G\ge K_{k}$.)  Conjecture~\ref{dcminor} remains open for all $t \ge 10$. It seems hard to generalize Theorem~\ref{dc9c}.

We need  some  known results  to prove our main results.  Before doing so, we need to define $(H,k)$-cockade. 
 For  a graph $H$ and an integer
$k$, let us define an {\it $(H,k)$-cockade} recursively as follows. Any graph isomorphic to $H$
is an $(H,k)$-cockade. Now let $G_1$, $G_2$ be $(H,k)$-cockades and let $G$
be obtained from the disjoint union of $G_1$ and $G_2$ by identifying a clique
of size $k$ in $G_1$ with a clique of the same size in $G_2$. Then the graph
$G$ is also an $(H,k)$-cockade, and every $(H,k)$-cockade can be constructed
this way.  We are now ready to state some known results.  The following  theorem is a result of Dirac~\cite{Dirac1964} for $p \le 5$ and Mader~\cite{mader} for $p \in \{6, 7\}$.
\begin{thm}\label{mader} (Dirac~\cite{Dirac1964}; Mader~\cite{mader})
For every integer $p\in \{1, 2, \dots, 7\}$, a graph on $n\ge p$ vertices
and at least $(p-2)n-{p-1\choose2}+1$ edges has a $K_p$ minor.
\end{thm}

J\o rgensen~\cite{Jorgensen1994} and later  Song and Thomas~\cite{SongThomas2006} generalized Theorem~\ref{mader} to $p=8$ and $p=9$, respectively, as follows.  

\begin{thm}\label{k8} (J\o rgensen~\cite{Jorgensen1994})
Every graph on $n \ge 8$ vertices with at least $6n - 20$ edges either contains a $K_8$-minor or is isomorphic to a $(K_{2,2,2,2,2}, 5)$-cockade.
\end{thm}

\begin{thm}\label{k9} (Song and Thomas~\cite{SongThomas2006})
Every graph on $n \ge 9$ vertices with at least $7n - 27$ edges either contains a $K_9$-minor, or is isomorphic to $K_{2,2,2,3,3}$, or is isomorphic to a $(K_{1,2,2,2,2,2}, 6)$-cockade.
\end{thm}

%It seems hard to generalize Theorem~\ref{mader} for all values of $p$.   In 2003, Seymour and Thomas~\cite{SongThomas2006} proposed the following conjecture. 
%
%\begin{conj}\label{conj2} (Seymour and Thomas~\cite{SongThomas2006})
%For every integer $p\ge1$ there exists a constant $N=N(p)$ such that
%every $(p-2)$-connected graph on $n\ge N$ vertices and at least
%$(p-2)n-{p-1\choose2}+1$ edges has a $K_p$ minor.
%\end{conj} 
%
%By Theorem~\ref{k9}, Conjecture~\ref{conj2} is true for $p\le9$.
%Norin and Thomas~\cite{norin} announced a proof of Conjecture~\ref{conj2} {\bf(to be updated)}. 

In our proof of Theorem~\ref{strongk89}, we need to examine graphs $G$ such that 
$k+1\le |G|\le 2k-5$, $\delta (G)\ge k-2$ and $G\not\ge K_k\cup K_1$. We shall use the following results. Lemma~\ref{k6k1} is a result   of J\o rgensen~\cite{Jorgensen1994}. 

\begin{lem}\label{k6k1} (J\o rgensen~\cite{Jorgensen1994})
Let $G$ be a graph with $n\le11$ vertices and  
$\delta(G)\ge6$ such that for every vertex $x$ in $G$, 
$G-x$ is not contractible to $K_6$. Then $G$ is one of the graphs $K_{2,2,2,2}, K_{3,3,3}$ or the complement of the 
Petersen graph.
\end{lem} 

Lemma~\ref{k6k1} implies  Lemma~\ref{k5k1} below. To see that, let $G$ be a graph  satisfying the conditions given  in Lemma~\ref{k5k1}.   By    applying Lemma~\ref{k6k1} to the graph obtained from  $G$  by adding $6-t$ vertices, each adjacent to every other vertex of the graph, we see that  $G\ge K_{t}\cup K_1$. 

\begin{lem}\label{k5k1}
For $t \in \{1, 2, 3, 4, 5\}$, let $G$ be a graph with $n\le 2t-1$ vertices and  
$\delta(G)\ge t$. Then $G\ge K_{t}\cup K_1$.  \end{lem} 

Lemma~\ref{k7k1} is  a result of  Song and Thomas~\cite{SongThomas2006}. Note that the proof of Lemma~\ref{k7k1} is computer-assisted. 

\begin{lem}\label{k7k1} (Song and Thomas~\cite{SongThomas2006})
Let $G$ be a graph with $|G| \in \{9, 10, 11, 12, 13\}$ such that $\delta(G)\ge7$. 
Then either  $G\ge K_7\cup K_1$,
%or $G$ is isomorphic to  $K_{1,2,2,2,2}$,
 or $G$ satisfies the following\medskip
 \myitem{(A)} 
either $G$ is isomorphic to $K_{1,2,2,2,2}$, or 
$G$ has four distinct vertices $a_1,b_1,a_2,b_2$ such that
 $a_1a_2, b_1b_2\notin E(G)$ and for $i=1,2$ 
the vertex $a_i$ is adjacent to $b_i$, the vertices $a_i, b_i$ have at most four common neighbors, and 
 $G+a_1a_2+b_1b_2\ge K_8$, 
\myitem{(B)} for any two sets  $A, B\subseteq V(G)$ of cardinality
 at least five such
 that neither is  complete and $A\cup B$ includes all vertices of $G$ of
degree at most $|G|-2$,
either
%\begin{itemize}
%\itemsep=-5pt
\myitemitem{(B1)}
 there exist $a\in A$ and $b\in B$ such that  $G'\ge K_8$, where $G'$ is
obtained from $G$ by adding all edges $aa'$ and $bb'$ for $a'\in A-\{a\}$ and
$b'\in B-\{b\}$, or
\myitemitem{(B2)}
 there exist $a\in A-B$ and $b\in B-A$ such that $ab\in E(G)$ and the
vertices $a$ and $b$ have at most five common neighbors in G,  or
\myitemitem{(B3)}
 one of $A$ and $B$ contains the other and $G+ab\ge K_7\cup K_1$ for all
distinct nonadjacent vertices $a,b\in A\cap B$.\\
% {\rm(A)} and {\rm(B)}.
\end{lem}

\section{Basic properties of non-complete double-critical graphs}

We begin with basic properties of non-complete double-critical $k$-chromatic graphs established  in~\cite{KPT2010}. We only list those that will be used in our proofs.\medskip

\begin{prop}\label{prop} (Kawarabayashi, Pedersen and Toft~\cite{KPT2010})
If $G$ is a non-complete double-critical $k$-chromatic graph, then the following hold:
\begin{enumerate}[(a)] \itemsep0pt \parskip0pt \parsep0pt \setlength{\itemindent}{1pt}

%\item $G$ does not contain a $K_{k-1}$-subgraph. \label{Kk-1}
\item $\delta(G) \ge k + 1$. \label{MinDeg}
%\item  For any $x\in V(G)$, $\alpha(N(x))%%\le d(x)-|B(xy)|-1
%\le d(x)-k+1$%, where $y\in N(x)$ is any vertex contained in an maximum independent set in $N[x]$.   \label{alpha}

%\item If $H$ is a connected subgraph of $G$, then the graph $G / H$ obtained by contracting $H$ to a single vertex is $(k-1)$-colorable. \label{Contraction}

\item Every edge $xy \in E(G)$ belongs to at least $k-2$ triangles. \label{Triangles}

%\item Every vertex $x \in V(G)$ has a neighbor $y$ such that $y$ is not complete to $N(x)$.  \label{A(xy)}

%\item There exists at least one edge $xy \in E(G)$ which is not a dominating edge of $G$. \label{D(xy)}

%\item If $x \in V(G)$ has neighbors $y, z$ say, such that $yz \notin E(G)$, then $x$ has another neighbor $w$ say, such that $wz \in E(G)$ and $wy \notin E(G)$.  \label{A(xy)edge}

%\item Any vertex $x$ with a non-neighbor in $G$ satisfies $\chi(N(x)) \le k- 3$. \label{cnbr}

%\item If $x \in V(G)$ is a vertex of degree $k+1$, then  $\overline{N(x)}$ consists only of isolated vertices and disjoint cycles of length at least five.  Moreover, $\overline{N(x)}$ contains at least one such cycle.  \label{complement}

%\item If $x,y \in V(G)$ are both of degree $k+1$, then $xy \notin E(G)$.  \label{min}
\item $G$ is 6-connected and no minimal separating set of $G$ can be partitioned into  two sets $A$ and $B$ such that $G[A]$ and $G[B]$ are edge-empty and complete, respectively.  \label{6-con}
\end{enumerate}
\end{prop}

%\begin{lem}\label{6-con}
%Suppose  $G$ is a non-complete double-critical $k$-chromatic graph.  Then $G$ is 6-connected and no minimal separating set of $G$ can be partitioned into two sets $A$ and $B$ such that $G[A]$ and $G[B]$ are edge-empty and complete, respectively.
%\end{lem}

Two proper vertex-colorings $c_1$ and $c_2$ of a graph $G$ are \dfn{equivalent}  if, for all $x, y\in V(G)$, $c_1(x)=c_1(y)$ iff $c_2(x)=c_2(y)$.   Two vertex-colorings $c_1$ and $c_2$ of a graph $G$ are \dfn{equivalent on a set  $A \subseteq V(G)$} if the restrictions ${c_1}_{|A}$ and ${c_2}_{|A}$ to $A$ are equivalent on the subgraph $G[A]$.  Let $S$ be a separating set of $G$, and let $G_1, G_2$ be connected subgraphs of $G$ such that $G_1 \cup G_2 = G$ and $G_1\cap G_2 = G[S]$.  If $c_1$ is a $k$-coloring of $G_1$ and $c_2$ is a $k$-coloring of $G_2$ such that $c_1$ and $c_2$ are equivalent on $S$, then it is clear that $c_1$ and $c_2$ can be combined to a $k$-coloring of $G$ by a suitable permutation of the color classes of, say $c_1$. The main technique in the proof of Proposition~\ref{prop}(\ref{6-con}) involves reassigning and permuting the colors on a separating set $S$ of a non-complete double-critical $k$-chromatic graph $G$ so that $c_1$ and $c_2$ are equivalent on $S$ to obtain a contradiction, where $c_1$ is a $(k-1)$-coloring of $G_1$ and $c_2$ is a $(k-1)$-coloring of $G_2$.  
It seems hard to use this idea to prove that every non-complete double-critical $k$-chromatic graph is $7$-connected, but we can use it to say a bit more about minimal separating sets of size $6$ in non-complete double-critical graphs.   
  %if there exists some permutation of the colors of $\phi_2$ such that $\phi_1(v) = \phi_2(v)$ for all $v \in V(G)$, i.e. the two partitions of $G$ induced by the color classes of $c_1$ and $c_2$ are equivalent up to the order of the partite sets.  

\begin{lem}\label{6-set}
Suppose  $G$ is a non-complete double-critical $k$-chromatic graph.  If $S$ is a minimal separating set of $G$ with $|S| = 6$, then either $G[S]\subseteq K_{3,\, 3}$  or $G[S]\subseteq K_{2,\, 2,\, 2}$.

\end{lem}

\pf
%Suppose $G$ is a non-complete double-critical $k$-chromatic graph.
By Propostion~\ref{prop}(\ref{6-con}), $G$ is 6-connected.  Let  $S = \{v_1, \dots, v_6\}\subset V(G)$ be a minimal separating set of $G$ such that neither $G[S]\subseteq K_{3,\, 3}$  nor $G[S]\subseteq K_{2,\, 2,\, 2}$.
Let $G_1$ and $G_2$ be subgraphs of $G$ such that $G_1 \cup G_2 = G$, $G_1 \cap G_2 = S$, and there are no edges from $G_1 - S$ to $G_2 - S$.
%Let $H$ be a component of $G - S$, and let $G_1=G[V(H)\cup S]$ and $G_2=G - V(H)$. Then $G_1 \cup G_2 = G$ and $G_1\cap G_2= S$.
Since $k \ge 6$ by Theorem~\ref{dc5c}, we have $\delta(G) \ge 7$ by Propostion~\ref{prop}(\ref{MinDeg}).
In particular, since $|S| = 6$, there must exist at least one edge $y_i z_i$ in $G_i - S$ for $i \in \{1, 2\}.$
It follows then that $G_i$ is $(k-2)$-colorable since it is a subgraph of $G - \{y_{3 - i} ,  z_{3 - i}\}$.
Let $c_1, c_2$  be $(k-2)$-colorings of $G_1 $ and $G_2$, respectively.  For $i=1,2$, define $|c_i(A)|$ to be the number of distinct colors assigned to the vertices of $A$ by $c_i$ for any $A\subseteq S$.
Clearly $c_1$ and $c_2$ are not equivalent on $S$, otherwise   $c_1$ and $c_2$, after a suitable permutation of the colors of $c_2$,  can be combined to  a $(k-2)$-coloring of $G$, a contradiction.
By Proposition~\ref{prop}(\ref{6-con}), $\alpha(G[S])\le4$ and so neither  $c_1$ nor $c_2$ applies the same color to more than four vertices of $S$.  Utilizing a new   color, say  $\beta$,  we next redefine the colorings $c_1$ and $c_2$  so that  $c_1$ and $c_2$ are $(k-1)$-colorings  of $G_1 $ and $G_2$, respectively,  and are equivalent on $S$. This yields a contradiction, as   $c_1$ and $c_2$, after a suitable permutation of the colors of $c_2$,  can be combined to  a $(k-1)$-coloring of $G$.  \medskip

Suppose that one of  the colorings $c_1$ and $c_2$, say $c_1$, assigns the same color to four vertices of $S$, say $c_1(v_3) = c_1(v_4) = c_1(v_5) = c_1(v_6)$. Then $\{v_3, v_4, v_5, v_6\}$ is an independent set in $G$.
By Proposition~\ref{prop}(\ref{6-con}), we must have $v_1 v_2 \notin E(G)$.
But then $G[S] \subseteq K_{2,\, 2,\, 2}$, a contradiction.
Thus neither $c_1$ nor $c_2$ assigns the same color to four distinct vertices of $S$. \medskip

Next suppose that one of  the colorings $c_1$ and $c_2$, say $c_1$, assigns the same color to three vertices of $S$, say $c_1(v_4) = c_1(v_5) = c_1(v_6)$. Then $\{v_4, v_5, v_6\}$ is an independent set in $G$.  Since $G[S]\not\subseteq K_{3,3}$, we have %If $c_2(v_1) = c_2(v_2) = c_2(v_3)$, then $\{v_1, v_2, v_3\}$ is an independent set. We redefine $c_2(v_4) = c_2(v_5) = c_2(v_6) = \beta$ and $c_1(v_1) = c_1(v_2) =c_1(v_3) = \beta$. Then $c_1$ and $c_2$ are equivalent on $S$ using $k-1$ colors. Thus  
 $|c_2(\{v_1, v_2, v_3\})| \ge2$. If  $|c_2(\{v_1, v_2, v_3\})| =2$, we may assume that $c_2(v_2) = c_2(v_3)$.
Then $\{v_2, v_3\}$ is an independent set. Then redefining $c_2(v_4) = c_2(v_5) = c_2(v_6) = \beta$ and $c_1(v_2) = c_1(v_3) = \beta$ will
%, after a suitable permutation of the colors of $c_2$,  
make $c_1$ and $c_2$ equivalent on $S$, a contradiction.
% using $k-1$ colors.
Thus $|c_2(\{v_1, v_2, v_3\})| =3$  and so $c_2$ assigns distinct colors to each of $v_1, v_2, v_3$.
We redefine  $c_2(v_4) = c_2(v_5) =c_2(v_6) = \beta$.
Clearly $c_1$ and $c_2$ are  equivalent on $S$  if  $c_1$ assigns distinct colors to each of $v_1, v_2, v_3$.
Thus $|c_1(\{v_1, v_2, v_3\})| \le2$.
Since $G[S]\not \subseteq K_{3,\, 3}$, we have $|c_1(\{v_1, v_2, v_3\})|=2$.
We may assume that  $c_1(v_2) = c_1(v_3)$.
Now  redefining  $c_1(v_3) = \beta$ yields that
%, after a suitable permutation of the colors of $c_2$,
$c_1$ and $c_2$  are equivalent  on $S$.
This proves that neither $c_1$ nor $c_2$ assigns the same color to three distinct vertices of $S$.
Thus $6\ge |c_i(S)| \ge 3$ $(i = 1, 2)$.
Since $G[S] \not\subseteq K_{2,\, 2,\, 2}$, we have $|c_i(S)| \ge 4$ $(i = 1, 2)$.
We may assume that  $|c_1(S)| \ge |c_2(S)|$.  
Then $|c_2(S)| \le 5$, for otherwise $c_1$ and $c_2$ are equivalent on $S$.
Thus  $ 5\ge|c_2(S)|\ge4$.    \medskip
 
Suppose that  $ |c_2(S)|=5$.  Then $ |c_1(S)|=5$ or $ |c_1(S)|=6$. We can make $c_1$ and $c_2$   equivalent on $S$ by assigning color $\beta$ to one of the two vertices that are colored the same color by $c_1$ (if $ |c_1(S)|=5$) and $c_2$.
Thus $ |c_2(S)|=4$.    Since neither $c_1$ nor $c_2$ assigns the same color to more than two distinct vertices of $S$, we may assume that $c_2(v_3)=c_2(v_4)$ and  $c_2(v_5) = c_2(v_6)$.
Then  $v_3v_4\notin E(G)$ and $ v_5v_6\notin E(G)$.
Since  $G[S]\not\subseteq K_{2,\, 2,\, 2}$, we have $v_1v_2\in E(G)$. Thus $c_1(v_1) \ne c_1(v_2)$.
We may assume that $c_1(v_3)\ne c_1(v_4)$ as $c_1$ and $c_2$ are not equivalent on $S$.
If $|c_1(S)| = 6$, then redefining $c_1(v_5) = c_1(v_6) = \beta$ and $c_2(v_3) = \beta$ will make $c_1$ and $c_2$ equivalent.
If $|c_1(S)| = 5$, then at least one of $v_3, v_4, v_5, v_6$ shares a color with another vertex of $S$, say $c_1(v_6) = c_1(v_i)$ for some $i \in \{1, \dots, 5\}$.
Then redefining $c_1(v_5) = c_1(v_6) = \beta$ and $c_2(v_3) = \beta$ will again make $c_1$ and $c_2$ equivalent.
Thus $|c_1(S)|=4$.
Suppose that one of $v_1$ or $v_2$ shares a color with another vertex of $S$.
Since $v_1 v_2 \in E(G)$, we may assume by symmetry that $c_1(v_1) = c_1(v_3)$.
If $c_1(v_5)$ and $c_1(v_6)$ are the two colors each assigned to only a single vertex of $S$ by $c_1$, then we also have $c_1(v_2) = c_2(v_4)$.
Now redefining $c_1(v_3) = c_1(v_4) = \beta$ and $c_2(v_5) = \beta$ will make $c_1$ and $c_2$ equivalent.
Hence one of the colors $c_1(v_5)$ and $c_1(v_6)$ is assigned to two vertices of $S$, say $c_1(v_6) = c_1(v_i)$ for some $i \in \{2, 4, 5\}$.
If $i = 2$ then redefine $c_1(v_5) = c_1(v_6) = \beta$ and $c_2(v_1) = c_2(v_3) = \beta$, if $i = 4$ then redefine $c_1(v_3) = c_1(v_4) = \beta$ and $c_2(v_6) = \beta$, and if $i = 5$ then redefine $c_1(v_3) = \beta$ and $c_2(v_3) = \beta$, and in each case $c_1$ is equivalent to $c_2$.
Therefore $c_1(v_1)$ and $c_1(v_2)$ are the two colors assigned to only a single vertex of $S$ by $c_1$.
Since $c_1$ and $c_2$ are not equivalent, we must have, say $c_1(v_3) = c_1(v_5)$ and $c_1(v_4) = c_1(v_6)$.
Now redefining $c_1(v_5) = c_1(v_6) = \beta$ and $c_2(v_3) = \beta$ will make $c_1$ and $c_2$ equivalent.
%
%
%Suppose $c_1(v_5)=c_1(v_6)$.
%Since $v_1v_2\in E(G)$ and $c_1$ and $c_2$ are not equivalent, we may assume that $c_1(v_1)=c_1(v_3)$.
%Now redefining  $c_1(v_3)=\beta$ will make $c_1$ and $c_2$ equivalent on $S$.
%Thus $c_1(v_5)\ne c_1(v_6)$.
%Let $A$ and $B$ be the two color classes of $c_1$ on $S$ with $|A|=|B|=2$.
%Suppose $v_1\in A$ and $v_2\in B$.
%We may assume that $v_3\in A$.
%Then $v_4\notin B$ because $G[S]\not\subseteq K_{2,\, 2,\, 2}$ and $v_1v_3\notin E(G)$.
%We may assume that $v_5\in B$. Redefining $c_1(v_5)=c_1(v_6)=\beta$ and $c_2(v_1)=\beta$ will make $c_1$ and $c_2$ equivalent on $S$.
%Thus $v_1, v_2\notin A\cup B$.
%By symmetry, we may assume $B=\{v_3, v_5\}$. We may further assume that $v_6\in A$.
%Now redefining $c_1(v_5)=c_1(v_6)=\beta$ will make $c_1$ and $c_2$ equivalent on $S$. 
%\medskip
%This completes the proof of Lemma~\ref{6-set}.   
\proofsquare

%\bsni Following the same process in an attempt to extend Theorem~\ref{6-set} to separating sets of cardinality 7 will result in the following.

%\begin{lem}\label{7-set}
%Suppose $G$ is a non-complete double-critical $k$-chromatic graph.  If $S$ is a minimal separating set of $G$ with $|S| = 7$, then either $G[S]\subseteq K_{3,\, 4}$  or $G[S]\subseteq K_{2,\, 2,\, 3}$  or $G[S]\subseteq K_{1,\,3,\, 3}$  or $G[S]\subseteq K_{1, \,2,\,2,\, 3}$. 
%\end{lem}

%\bsni The proof of Lemma~\ref{7-set} is similar to the proof of Lemma~\ref{6-set}, and is omitted.  Additionally, the presence of an independent set of cardinality 1 in conclusions (iii) and (iv) of Lemma~\ref{7-set} allow for a vertex $v$ in our separating set $S$ with no required adjacency or non-adjacency relations in $S$.  In particular, $v$ is potentially complete to $S \setminus \{v\}$.  In practice, this will mean that there is no benefit in utilizing Lemma~\ref{7-set} over Lemma~\ref{6-set}.

%&&&&&&&&&&&&&&&&&&
%&&&&&&&&&&&&&&&&&
\section{Proofs of Theorem~\ref{strongk89} and Theorem~\ref{dc9c}}

 In this section we first prove Theorem~\ref{strongk89}.\medskip
 
\pf  Let $G$ be  a graph as in the statement with $n$ vertices.    By assumption, we have 
\medskip

\noindent(1) 
$k+1\le \delta(G)\le 2k-5$ and $\delta(N(x))\ge k-2$ for any $x$ in $G$; and \medskip

\noindent(2)  $G$ is $(k-3)$-connected and for any   minimal separating set $S$ of $G$ and any $x\in S$, $G[S\less \{x\}]$   is not  a complete subgraph. \medskip

We first show that the statement is true for $k=6$.
Then $G$ is $3$-connected with $\delta(G)=7$.  The statement is trivially true if $G$ is complete, so we may assume $G$ is not complete.
Let $x\in V(G)$ be a vertex of degree 7.
By (1),  $\delta(N(x))\ge 4$, and so $e(N(x))\ge14$.  If $e(N(x))\ge16$, then by Theorem~\ref{mader}, $N(x)\ge K_5$ and so $G\ge N[x]\ge K_6$.  If $e(N(x)) = 15$, then let $K$ be a component of $G-N[x]$ with $|N(K)|$  minimum. 
By (2),  $|N(K)|\ge3$ and $N(K)$ is not complete.  Let $y,z\in N(K)$ be non-adjacent in $N(x)$ and let $P$ be a $(y,z)$-path with interior vertices in $K$. We see that $G\ge K_6$ by contracting all but one of the edges of $P$. 
So we may assume that $e(N(x))=14$, and so $N(x)$ is $4$-regular and  $\overline{N(x)}$ is 2-regular.  Thus  $\overline{N(x)}$ is then either isomorphic to $C_7$ or to $C_4 \cup C_3$, and in both cases it is easy to see that $N(x) \ge K_5$ and thus $G \ge K_6$, as desired.  Hence we may assume $7\le k\le9$.

%We first show that the statement is true for $k=6$. Assume $k=6$.  Then $G$ is $3$-connected with $\delta(G)=7$.  Let $x\in V(G)$ be  a vertex of degree $7$.   By (1),  $\delta(N(x))\ge 4$, and so $e(N(x))\ge14$.  Assume that $e(N(x))\ge15$. By Theorem~\ref{mader}, $N(x)>K_5$ and so $G>N[x]>K_6$. So we may assume that $e(N(x))=14$,  and so $N(x)$ is $4$-regular and $G-N[x]$ is non-null. Let $K$ be a component of $G-N[x]$.   By (2),  $|N(K)|\ge3$ and $N(K)$ is not complete.  Let $x,y\in N(K)$ be non-adjacent in $N(x)$ and let $P$ be an $(x,y)$-path with interior vertices in $K$. We see that $G>K_6$ by contracting all but one of the edges of $P$, as desired. \medskip  So we may assume $7\le k\le8$. 

Suppose for a contradiction that $G\not\ge K_k$. We next prove the following.\medskip
 
\noindent  (3)  Let $x\in V(G)$ be  such that $k+1\le d(x)\le 2k-5$. Then there is no
 component $K$ of  $G-N[x]$ such that $N(K')\cap M\subseteq N(K)$ for every component $K'$
 of $G-N[x]$, where $M$ is the set of  vertices of $N(x)$ 
not adjacent to  all  other vertices 
of $N(x)$.\medskip

\pf  Suppose   such a component $K$ exists. Among all vertices  $x$ 
with $k+1\le d(x)\le 2k-5$ for which such a component exists, choose $x$ to be of minimal degree, and among all such components $K$ of $G - N[x]$, choose $K$ such that $|N(K)|$ is minimum.
We first prove 
that $M\subseteq N(K)$. 
Suppose for a contradiction
that $M- N(K)\ne\emptyset$,
 and let $y\in M\setminus N(K)$ be such that $d(y)$ is minimum. Clearly, $d(y)< d(x)$.  
Let $J$ be the component of $G-N[y]$ containing $K$.  Since $d(y)<d(x)$
the choice of $x$ implies that  $N(x)\setminus N[y]\not\subseteq V(J)$. Let  $H=N(x)\setminus (N[y] \cup N(K))$. 
We have $d_G(z)\ge d_G(y)$ for all $z\in V(H)$ by the choice of $y$.
Let $t=|V(H)|$.  Then $t\ge 2$, for otherwise the vertex $y$ and
component $H$ contradict the choice of $x$.  On the other hand
$t\le d(x)-d(y)\le (2k-5)-(k+1)=k-6\le3$ and so $k\ge8$. Notice that $t=2$ when $k=8$. From (1) applied to $y$
we deduce that $N(y)\cap N(x)$ has minimum degree at least $k-3$.
Let $L$ be the subgraph of $G$ induced by $(N[y]\cap N(x))\cup V(H)$.
Then the edge-set of $L$ consists of edges of $N(x)\cap N(y)$, edges incident
with $y$, and edges incident with $V(H)$.  Clearly, $e(L - V(H), H)=\sum_{z\in V(H)} \, (d(z)-1) -2e(H)\ge t(d(y)-1)-2e(H)$. Thus
\begin{align*}
e(L)&\ge \frac{(k-3)\, (d(y)-1)}2+d(y)-1+e(L - V(H), H)+e(H)\\
&\ge \frac{(k-3)\, (d(y)-1)}2+d(y)-1+t(d(y)-1)-e(H)\\
&\ge \frac{(k-3)(d(y)-1)}2+d(y)-1+t(d(y)-1)-\frac12\, t(t-1)\\
& \ge \left\{
       \begin{array}{ll}
           5(d(y)+2)+\frac{d(y)}2-\frac{33}{2} & \quad  {\text{if}} \quad k=8 \\
           6(d(y)+t)+(t-2)d(y)-4-7t-\frac12 \, t(t-1) & \quad  {\text{if}} \quad k=9
      \end{array}
    \right.\\
&   \ge (k-3)|V(L)|-{k-2\choose2}+1,
\end{align*}
because $d(y)\ge k+1$  and $2\le t\le k-6$. If $k=9$, since $12\le |V(L)|\le 13$ the 
graph $L$ is not a $(K_{2,2,2,2,2},5)$-cockade.  By Theorem~\ref{mader} and Theorem~\ref{k8},  $N(x)\ge L\ge K_{k-1}$.  Thus $G\ge N[x]\ge K_k$, a contradiction.  This proves
that $M\subseteq N(K)$.\medskip

If  $N(x)\ge K_{k-2}\cup K_1$, 
then $N(x)$ has a vertex $y$ such that $N(x)-y\ge K_{k-2}$. If $y\not\in M$, then
$N(x)\ge K_{k-1}$. Otherwise,
by contracting the connected set $V(K)\cup \{y\}$ we can contract $K_{k-1}$ onto
$N(x)$.  Thus in either case
$G\ge K_{k}$, a contradiction. 
Thus $N(x)\not\ge K_{k-2}\cup K_1$.  If $k\le 8$, by  Lemma~\ref{k6k1} and Lemma~\ref{k5k1},  we have $k=8$ and $N(x)$ is  either $K_{3,3,3}$ or $\overline{P}$, where $\overline{P}$ is the complement of the Petersen graph. If  $N(x)=\overline{P}$, it can be easily checked that $\overline{P}+yz \ge K_7$ for any $yz \in E(P)$.
By (2),  $|N(K)|\ge5$ and $N(K)$ is not complete.
Let $y, z\in N(K)$ be non-adjacent in $N(x)$ and let $Q$ be a $(y,z)$-path with interior vertices in $K$.
We see that $G\ge K_8$ by contracting all but one of the edges of $Q$, a contradiction.
Thus  $N(x)= K_{3,3,3}$, and so $M = N(x)$.
Let $\{a_1, a_2,a_3\}$ and  $\{b_1, b_2,b_3\}$ %  and $\{c_1, c_2,c_3\}$ 
be the vertex sets of two disjoint triangles of $\overline{N(x)}$.  
 Suppose $G-N[x]$ is 2-connected or has at most two vertices.
By Proposition~\ref{prop}(\ref{Triangles}), the vertices $a_i, b_i$ (i=1,2) have at least two common neighbors in $G-N[x]$.  Let $u_1, u_2$  (resp. $w_1, w_2$) be two  distinct  common neighbors
of $a_1$ and $b_1$ (resp. $a_2$ and $b_2$) in $G-N[x]$.  By Menger's Theorem,
$G-N[x]$ contains two disjoint paths from $\{u_1, u_2\}$ to $\{w_1, w_2\}$ and
so $G \ge N[x]+a_1a_2+ b_1b_2 \ge K_8$, a contradiction. Thus  $G-N[x]$ has at least three vertices and is
not $2$-connected. If $G-N[x]$ is disconnected, let $H_1=K$ and $H_2$  be another  connected
component of $G-N[x]$. If $G-N[x]$ has a cut-vertex, say $w$,  let $H_1$  be a
connected component of $G-N[x]-w$  and let  $H_2=G-N[x]-V(H_1)$.
In either case, $H_1$ and $H_2$ 
are disjoint connected subgraphs of $G-N[x]$ such that $M\subseteq
N(H_1)\cup N(H_2)$ (because we have shown that $M\subseteq N(K)$). Thus  $N(H_1)\cup N(H_2)=N(x)$ because $M=N(x)$. 
By (2), $N(H_i) $ is not complete and $|N(H_i) | \ge
4$ since $k = 8$.
   Thus  each of $N(H_1)$ and $N(H_2)$ must contain at least one edge of $\overline{N(x)}$. 
%If for some $i \in \{1, 2\}$, $N(H_i)$ contains only one edge $e$ of $\overline{N(x)}$, then $|N(H_i)| = 4$ and $N(H_{3 - i})$ contains at least one edge of $\overline{N(x)}$ disjoint from $e$. If each of $N(H_1)$ and $N(H_2)$ contains at least two edges of $\overline{N(x)}$, then it is not hard to see that we can always select an edge of $\overline{N(x)}$ in $N(H_1)$ and an edge of $\overline{N(x)}$ in $N(H_2)$ such that the selected edges are disjoint.
%
%
%we may assume that  $\overline {N(H_2)}$ contains  a matching of size at least two.
Since $N(x)=K_{3,3,3}$ and $N(H_1)\cup N(H_2)=N(x)$, we may thus assume that $a_1a_2\in \overline {N(H_1)}$ and  $b_1 b_2\in\overline {N(H_2)}$. By contracting $H_1$ onto $a_1$ and $H_2$ onto $b_1$ we see that     $G\ge N[x]+a_1a_2+b_1b_2\ge K_8$, a contradiction. 
This proves that   $k=9$ and so  by Lemma~\ref{k7k1}, 
we may assume that $N(x)$ 
satisfies properties (A) and  (B). 
\medskip

Since $d(x)\ge10$, $N(x)\ne K_{1,2,2,2,2}$. If $G-N[x]$ is 2-connected or has at most two vertices, 
then by property (A) and (2),  the set $N(x)$ has four distinct
vertices  $a_1, b_1, a_2, b_2$ such that $a_1a_2, b_1b_2\notin E(G)$,
 $N(x)+a_1a_2+b_1b_2 \ge K_8$ and for 
$i=1,2$ the vertex $a_i$ is adjacent to $b_i$, and the vertices $a_i, b_i$  
have at least two common neighbors in
$G-N[x]$. Let $u_1, u_2$  (resp. $w_1, w_2$) be two  distinct  common neighbors
of $a_1$ and $b_1$ (resp. $a_2$ and $b_2$) in $G-N[x]$.  By Menger's Theorem,
$G-N[x]$ contains two disjoint paths from $\{u_1, u_2\}$ to $\{w_1, w_2\}$ and
so $G \ge N[x]+a_1a_2+ b_1b_2 \ge K_9$, a contradiction.
Thus  $G-N[x]$ has at least three vertices and is
not $2$-connected.
If $G-N[x]$ is disconnected, let $H_1=K$ and $H_2$  be another  connected
component of $G-N[x]$. If $G-N[x]$ has a cut-vertex, say $w$,  let $H_1$  be a
connected component of $G-N[x]-w$  and let  $H_2=G-N[x]-V(H_1)$.
In either case, $H_1$ and $H_2$ 
are disjoint connected subgraphs of $G-N[x]$ such that $M\subseteq
N(H_1)\cup N(H_2)$ (because we have shown that $M\subseteq N(K)$).
% and if $G-N[x]$ is connected,
%we have $V(H_1)\cup V(H_2)=V(G)-N[x]$.
For $i=1,2$ let $A_i=N(H_i)\cap N(x)$.
By (2), $A_i$ is not complete and $|A_i|\ge
5$ for $i=1,2$.
%We will only use these properties of $H_1$ and $H_2$, and
%therefore may use the symmetry between $H_1$ and $H_2$.
By property (B),
$A_1$ and $A_2$ satisfy properties (B1), (B2) or (B3).\medskip

 Suppose first that $A_1$ and $A_2$ satisfy property (B1). Then there exist
$a_i\in A_i$  such that   $N(x)+\{a_1a: a\in
A_1 \setminus \{a_1\}\}+\{a_2a: a\in A_2 \setminus \{a_2\}\}\ge K_8$. By contracting the
connected sets
$V(H_1)\cup \{a_1\}$ and $V(H_2)\cup\{a_2\}$ to single  vertices,  we
see that  $G\ge K_9$, a contradiction.
Suppose next that $A_1$ and $A_2$ satisfy
property (B2).  Then there exist $a_1\in A_1\setminus A_2$ and $a_2\in A_2\setminus A_1$
such that $a_1a_2\in E(G)$ and  the vertices  $a_1$ and $a_2$ have  at
most five common neighbors in $N(x)$. Thus $a_1,a_2\in M$  by
(1), and by another application of (1)
there exists a common neighbor $u\in V(G)\setminus N[x]$
of $a_1$ and $a_2$. But $a_1\not\in A_2$ and $a_2\not\in A_1$,
and hence $u\not\in V(H_1)\cup V(H_2)$. Thus $G-N[x]$ is disconnected
and $H_1=K$. But then $a_2\in M\subseteq N(K)=N(H_1)$, a contradiction.
Thus we may assume that $A_1$ and $A_2$ satisfy (B3), and hence
$A_i\subseteq A_{3-i}$ for some $i\in\{1,2\}$.
As $M\subseteq A_1\cup A_2$, we have $M\subseteq N(H_{3-i})$.
Since $A_i$ is not complete, let $a,b\in A_i$ be distinct and  not adjacent.
By property (B3),  $N(x)+ab \ge K_7\cup K_1$. Let
$P$ be an $(a, b)$-path with interior in $H_i$.
By contracting all but one of the
edges of the path $P$ and by contracting $H_{3-i}$ similarly as above,
we see that $G\ge K_9$, a contradiction.
\hfill\vrule height3pt width6pt depth2pt \\

\noindent  (4)  $G-N[x]$  is disconnected for
every vertex $x\in V(G)$ of degree at most $2k-5$.

\pf If $G-N[x]$ is not null, then it is disconnected by (3). 
Thus we may assume that $x$ is adjacent to every other vertex of $G$. Let $H=G-x$. 
 Then $|H|=d(x)$ and $\delta(H)\ge k$. Thus $e(H)\ge \frac{k\, d(x)}2> (k-3)\, d(x)-{k-2\choose2}+1$ because $d(x)\le 2k-5$. By Theorem~\ref{mader} and Theorem~\ref{k8},  
$G-x$ has a $K_{k-1}$ minor and so  the graph $G$ has a $K_k$ minor, a contradiction.
\hfill\vrule height3pt width6pt depth2pt \\

\noindent {(5)} Let $x\in V(G)$ be  such that $k+1\le 
d(x)\le 2k-5$. 
Then there is no component $K$ of $G-N[x]$ such that $d_G(y)\ge2k-4$ for every vertex 
$y\in V(K)$.

\pf Assume that such a component $K$ exists. 
%By Lemma~\ref{lem35}, $|K|\ge3$. 
 Let 
$G_1=G-V(K)$ and $G_2=G[V(K) \cup N(K)]$. Let $d_1$ be the  maximum number of edges that can be added to $G_{2}$ by contracting edges of $G$ with at least one 
end in $G_1$. More precisely, let $d_1$ be the 
largest integer so that $G_1$ contains disjoint sets of vertices 
$V_1, V_2, \dots, V_p$ so that $G_1[V_j]$ is connected, 
 $|N(K)\cap V_j|=1$ for  $1\le j\le p =|N(K)|$, and so that the graph obtained from $G_1$ by contracting $V_1, V_2, \dots, V_p$ 
and deleting $V(G)\setminus (\bigcup_j V_j)$ has 
$e(N(K))+d_1$ edges.
  Let $G_2'$
 be a graph with $V(G_2')=V(G_2)$ and  $e(G_2')=e(G_2)+d_1$ edges 
obtained from $G$ by contracting edges in $G_1$. 
By (1), $|G_2'|\ge k+2$.  If
 $e(G_2')\ge (k-2)\, |G_2'|-{k-1\choose2}+2$, then by Theorem~\ref{mader} and Theorem~\ref{k8},  $G\ge G_2' \ge  K_k$, 
a contradiction.
 Thus 
%\begin{equation}
$$
e(G_2)=e(G_2')-d_1\le(k-2)\, |G_2|-{k-1\choose2}+1-d_1=
(k-2)|N(K)|+(k-2)|K|- {k-1\choose2}+1 -d_1.
$$
%\end{equation} 
By contracting the edge $xz$, 
where $z\in N(K)$ has minimum degree $d$ in $N(K)$, we see that $d_1\ge 
|N(K)|-d-1$ %, where $d=\delta(N(K))$, 
and hence
\begin{equation}\tag{a}
e(G_2)\le
(k-3)|N(K)|+(k-2)|K|- {k-1\choose2}+2+d.
\end{equation}

Let $t=e_G(N(K), K)$. We have 
$e(G_2)= e(K)+t+e(N(K))$ and 
\begin{equation}\tag{b}
2e(K)\ge (2k-4)|K|-t,
\end{equation}
and hence 
\begin{equation}\tag{c}
e(G_2)
\ge(k-2)|K|+t/2+ d|N(K)|/2.
\end{equation}

Since $N(x)$ has minimum 
degree at least $k-2$, it follows that 
the subgraph $N(K)$ of $N(x)$ has minimum degree at least $(k-2)-(d(x)-|N(K)|)$.
Thus $d\ge(k-2)-(d(x)-|N(K)|)\ge|N(K)|-k+3$.
 From (a) and (c) we get
\begin{equation}\tag{d}
-t/2\ge 
-(k-3)|N(K)|+d(|N(K)|-2)/2+{k-1\choose2}-2\ge  \left\{
       \begin{array}{ll}
       -8 & \quad {\text{if}} \quad k=7 \\
           -14& \quad  {\text{if}} \quad k=8\\
           -18 & \quad  {\text{if}} \quad k=9
      \end{array}
    \right.
\end{equation}
where  the second inequality becomes $\frac{t}2\le 11$ when $|N(K)|=2k-6$ and  $k=7,8$, and   the second inequality holds 
with equality only when $|N(K)|=10$ and $k=9$.
Since $G$ is not contractible to $K_k$,  we deduce from (b) and Theorem~\ref{mader}, Theorem~\ref{k8} and Theorem~\ref{k9} that 
$|K|<8$. The inequalities $e(K)\ge  5|K| -8 $ when $k=7$,   
       $e(K)\ge    6|K|-14$ when $k=8$, and $e(K)\ge7|K|-18$ when $k=9$
     imply $|K|\le3$. But
every vertex of $K$ has degree at least $2k-4$ and $N(K)$ is a proper
subgraph of $N(x)$, and hence $|K|=3$, $|N(K)|=2k-6$ and $\frac{t}2=3(k-3)\ge12$ when $k=7, 8$, and (d) holds with equality for $|N(K)|=12$ when $k=9$, contrary to our earlier observation of (d)  that $\frac{t}2\le 11$ when $|N(K)|=2k-6$ and $k=7,8$, and (d) holds with equality only when $|N(K)|=10$ and $k=9$. \hfill\vrule height3pt width6pt depth2pt\bigskip

By (1) 
there is a vertex $x$ of degree $k+1, k+2, \dots,$  or $2k-5$ in $G$.
Choose such a vertex $x$ so that $G-N[x]$ has a component $K$ of minimum order. 
Then choose a vertex $y\in V(K)$ of least degree in $G$. 
Thus $k+1\le d_G(y)\le 2k-5$ by (1) and (5). Let $L$ be the component of 
$G-N[y]$ containing $x$. We claim that  $N(L)$ contains all vertices of $N(y)$ that are not
adjacent to all other vertices of $N(y)$.
 Indeed,  let $z\in N(y)$ be not adjacent to some vertex of
$N(y)\setminus \{z\}$. We may assume that $z\notin N(x)$, for otherwise $z\in N(L)$.
Thus $z\in V(K)$, and hence $d_G(z)\ge d_G(y)$ by the choice of $y$. Thus
$z$ has a neighbor $z'\in N[x]\cup V(K) \setminus N[y]$.  Then $z'\in V(L)$, 
for otherwise the component of $G-N[y]$ containing $z'$  would be a proper
subgraph of $K$. Thus  $z\in N(L)$. 
This proves our claim that $N(L)$ contains all vertices $z$ as
above, contrary to (3). This contradiction completes the proof of
  Theorem~\ref{strongk89}.
\hfill\vrule height3pt width6pt depth2pt \\

We are now ready to prove Theorem~\ref{dc9c}.\medskip

\pf  Let $G$ be a double-critical $t$-chromatic graph with $t\ge k$. The assertion is trivially true if $G$ is complete. By Theorem~\ref{dc5c},  we may assume that $t\ge6$.  By Proposition~\ref{prop}(\ref{MinDeg}), $\delta(G)\ge k+1$. By Theorem~\ref{mader}, Theorem~\ref{k8} and Theorem~\ref{k9}, we have $\delta(G)\le 2k-5 $.
By Proposition~\ref{prop}(\ref{Triangles}),  every edge of $G$ is contained in at least $k-2$ triangles.   By Proposition~\ref{prop}(\ref{6-con}), $G$ is $6$-connected and no minimal separating set of $G$ can be partitioned into a  clique and an independent set. By Theorem~\ref{strongk89}, $G\ge K_k$, as desired. \hfill\vrule height3pt width6pt depth2pt

\section*{Acknowledgement}
The authors would like to thank the anonymous referees for many helpful comments.

%**************************
%&&&&&&&&&&&&&&&&

\end{document}